# Generic Simplicity of a Schrödinger-type Operator on the Torus

## Louis Omenyi[1*], Emmanuel Nwaeze[1] and McSylvester Omaba[1]

[1]*Department of Mathematics, Computer Science, Statistics and Informatics, Federal University, Ndufu-Alike, Ikwo, Nigeria.*

*Authors' contributions*

*This work was carried out in collaboration between all authors. Author LO designed the study, performed the mathematical analysis and wrote the first draft of the manuscript. Author EN managed the analyses and editing of the study. Author MO managed the literature searches. All authors read and approved the final manuscript.*

| *Original Research Article* |  |
|---|---|

## Abstract

The generic simplicity of the spectrum of a Schrödinger-type operator on the n-dimensional torus is studied using the Rayleigh-Schrödinger perturbation theory. The existence of a perturbation potential of the Laplacian is proved and suitable conditions on the potential that guarantee the generic simplicity of the spectrum constructed. It is also proved that with the potential, the degeneracy of the spectrum of the Laplacian on the n-dimensional torus splits at first order of the perturbation.

## 1 Introduction

Pierre-Simon Laplace, according to literature, e.g. [1], discovered that a gravitational field can be represented as the gradient $\nabla u = \dfrac{\partial u}{\partial x}$ of a potential function $u$. He further showed that in a free space, $\nabla \cdot \nabla u = \Delta u = 0$. This equation later became known as Laplace equation and $\Delta$ known as the Laplace operator or simply the "Laplacian", [2].

---

*\*Corresponding author: E-mail: omenyi.louis@funai.edu.ng;*

Closely related to the Laplacian is the Schrödinger operator. Consider for a general quantum system, the usual Schrödinger equation given by

$$i\hbar \frac{\partial}{\partial t}\psi = \hat{H}\psi$$

where $\psi$ is the wave function, $i\hbar \dfrac{\partial}{\partial t}$ the energy operator and $\hat{H}$ the usual Hamiltonain operator; see e.g. [3,4,5,6] among other literature. With a single particle perturbation, the equation turns out to be

$$i\hbar \frac{\partial}{\partial t}\psi(r,t) = -\frac{\hbar^2}{2m}\Delta\psi(r,t) + V(r)\psi(r,t)$$

where $-\dfrac{\hbar}{2m}\Delta$ is the kinetic energy operator, $\dfrac{\hbar}{2m}$ the Planck constant, $m$ the mass of the energy function and $V(r)$ the time-dependent potential energy at the position $r$ with $\psi(r,t)$ the probability amplitude for the particle to be found at position $r$ at time $t$; [7]. The operator $\hat{H}$ is given by

$$\hat{H} = -\frac{\hbar^2}{2m}\Delta + V$$

which is of the form usually called Schrödinger operator in literature, c.f [2,8] and [7].

In this paper, we simplify the coefficients of the Schrödinger operator $\hat{H}$ by reducing its coefficients to unity, in some sense. Thus we have the resulting linear operator which we denote by $H$ and written as

$$H = \Delta + V \tag{1}$$

where, again, $\Delta$ is the Laplacian and $V$ is a perturbation potential. We call $H$ a Schrödinger-type operator.

The main goal of studying perturbation theory of the Schrödinger-type operator in this paper is to attempt to specify the nature of the potential $V$ of (1) that preserves self-adjointness and splits the spectrum of $H$ on $n$-dimensional unit torus.

Let $(M, g)$ be a closed connected smooth Riemannian manifold. The Laplacian on $C^\infty(M)$ is the operator

$$\Delta_g : C^\infty(M) \to C^\infty(M) \tag{2}$$

defined in local coordinates by

$$\Delta_g = -\operatorname{div}(\operatorname{grad}) = -\frac{1}{\sqrt{|g|}}\sum_{i,j}\frac{\partial}{\partial x^i}(\sqrt{|g|}\,g^{ij}\frac{\partial}{\partial x^j}). \tag{3}$$

The operator $\Delta_g$ extends to a self-adjoint operator on $L^2(M) \supset H^2(M) \to L^2(M)$ with compact resolvent, [9,10] and [11]. This implies that there exists orthonormal basis $\{\psi_k\} \subset L^2(M)$ consisting of eigenfunctions such that

$$\Delta_g \psi_k = \lambda_k \psi_k \qquad (4)$$

where the eigenvalues $\{\lambda_k\}_{k=1}^{\infty}$ are listed with multiplicities; [12,9,11] and [10].

The question of generic simplicity of these eigenvalues on different Riemannian manifolds has been tackled extensively in Riemannian geometry and related fields using various perturbation techniques. For example, since the pioneering works of Rayleigh and Lindsay [13] and Hadamard [14], various kinds of perturbations have been carried-out in order to split the spectrum; see e.g. [15,16,17,18,5,19,6,20,21]; more recently [22,23] and still ongoing [24]. In this paper, we prove the existence of a perturbation potential $V$ of the Laplacian that guarantees the simplicity of the spectrum, at first order, on the $n$-dimensional unit torus. We proceed by fixing notations and basic concepts.

## 2 The $n$-dimensional Torus

One may think of a point on the unit circle, $S^1$, in the Cartesian plane, as the pair $(x, y)$ with $x^2 + y^2 = 1$. Similarly, a point in the 2-dimensional unit torus, $S^1 \times S^1$, is naturally described by co-ordinate pairs $(x_1, y_1), (x_2, y_2)$ such that $x_1^2 + y_1^2 = 1$ and $x_2^2 + y_2^2 = 1$; [12]. By extension, an n-dimensional flat torus (or simply n-torus) is the product of n-circles; [1]:

$$\mathsf{T}^n = \underbrace{S^1 \times S^1 \times \cdots \times S^1}_{\text{n--times}}.$$

In general, the $n$-torus is an Abelian Lie group, [12]. A point on $\mathsf{T}^n$ is represented by $n$-tuple $(\phi_1, \phi_2, \cdots, \phi_n)$ where each $\phi_k$ is understood as an "angle" defined modulo $2\pi$; [12]. The binary operation on this group is simply addition modulo $2\pi$, that is

$$(\phi_1, \phi_2, \cdots, \phi_n) + (\psi_1, \psi_2, \cdots, \psi_n) = ((\phi_1 + \psi_1), (\phi_2 + \psi_2), \cdots, (\phi_n + \psi_n)) \bmod 2\pi.$$

The $n$-dimensional torus $\mathsf{T}^n$ is fully represented as a matrix Lie group by simply assigning to each $n$-tuple $(\phi_1, \phi_2, \cdots, \phi_n)$ the following diagonal matrix [7,3,12]:

$$A_\phi = \begin{pmatrix} \mathrm{e}^{i\phi_1} & 0 & 0 & 0 \\ 0 & \mathrm{e}^{i\phi_2} & 0 & 0 \\ \vdots & \vdots & \vdots & \vdots \\ 0 & ... & ... & \mathrm{e}^{i\phi_n} \end{pmatrix}. \qquad (5)$$

# 3 Some Calculus on the $n$-dimensional Torus

From the foregoing, it is clear that $\mathsf{T}^n$ is a product Riemannian manifold. So, let $(X, g_X)$ and $(Y, g_Y)$ be smooth Riemannian manifolds, it is known that the product manifold $(M, g)$ with $M = X \times Y$ has Riemannian metric $g = g_X + g_Y$; [12,1,9,10]. In matrix form

$$g = \begin{pmatrix} (g_X) & 0 \\ 0 & (g_Y) \end{pmatrix}. \tag{6}$$

Consequently, the Lapacian on the $n$-torus with the standard Cartesian metric is

$$\Delta = \frac{\partial}{\partial x_1} + \frac{\partial}{\partial x_2} + \cdots + \frac{\partial}{\partial x_n}. \tag{7}$$

It is known, (see e.g. [12,1,9,10,25]), that the spectrum of the Laplacian on $C^\infty(\mathsf{T}^n)$ is given by

$$\sigma(\Delta) = k_1^2 + k_2^2 + \cdots + k_n^2; k \in \mathsf{Z}^n \tag{8}$$

with multiplicity $m_j < 2j + 4$ for the $j^{th}$ eigenvalue.

Integration can also be performed over $\mathsf{T}^n$. It is to be understood as [12,1]:

$$\int_{\mathsf{T}^n} f(x)dx = \int_{S^1}\int_{S^1}\cdots\int_{S^1} f(x_1, x_2, \cdots, x_n)dx_1 dx_2 \cdots dx_n. \tag{9}$$

That is,

$$\int_{\mathsf{T}^n} f(x)dx = \int_0^{2\pi}\int_0^{2\pi}\cdots\int_0^{2\pi} f(x_1, x_2, \cdots, x_n)dx_1 dx_2 \cdots dx_n \tag{10}$$

The space $L^p(\mathsf{T}^n)$, [9,1], consists of all Lebesgue measurable functions such that

$$\| f \|_{L^p(\mathsf{T}^n)}^p := \int_{\mathsf{T}^n} | f(x) |^p \, dx < \infty. \tag{11}$$

When $p = 2$, we say that the function $f$ is absolutely square integrable. For $f, g \in \mathsf{T}^n$, we define

$$\langle f, g \rangle = \int_{\mathsf{T}^n} \bar{f}(x) g(x) dx. \tag{12}$$

Given n-tuples of integers $k = (k_1, k_2, \cdots, k_n) \in \mathsf{Z}^n$, and $x \in \mathsf{T}^n$; we form the scalar

$$k.x = k_1 x_1 + k_2 x_2 + \cdots + k_n x_n = \sum_{j=1}^n k_j x_j. \tag{13}$$

For any $2\pi$ -periodic smooth function $f \in C^\infty(\mathsf{T}^n)$, [1], define (for $s \geq 0$ an integer), the space:

$$\mathsf{H}^s(\mathsf{T}^n) = \{f \in L^2(\mathsf{T}^n) : \sum_{k \in \mathsf{Z}^n} (1 + |k|^2)^s |\hat{f}(k)|^2 < \infty\}. \tag{14}$$

$\mathsf{H}^s(\mathsf{T}^n)$ is the Sobolev space on $\mathsf{T}^n$. Note, $\hat{f}$ is the Fourier transform of $f \in \mathsf{T}^n$; see e.g. [1].

**Theorem 3.1** [1,9]. *The complex-valued function* $e_k : \mathsf{T}^n \to \mathsf{C}$ *given by*

$$e_k(x) = (2\pi)^{-\frac{n}{2}} \prod_{j=1}^{n} e^{ik_j \cdot x_j}$$

*constitute an orthonormal basis in* $L^2(\mathsf{T}^n, \mathsf{C})$ *and hence an orthogonal basis in the Sobolev space* $\mathsf{H}^p(S^1)$.

We immediately have the following preparatory and standard result given as lemma (3.2) below.

**Lemma 3.2.** *The Laplacian* $\Delta$ *is self-adjoint on* $\mathsf{H}^2(\mathsf{T}^n)$.

*Proof.* The Laplacian is naturally defined in $\mathsf{H}^2(\mathsf{T}^n)$ in Fourier space as a multiplication operator [8]. Thus, for $f \in L^2(\mathsf{T}^n)$, we have

$$\hat{\Delta f}(k) = \|k\|^2 \, \hat{f}(k)$$

where of course $\|k\|^2 \in L^2(\mathsf{T}^n)$. For any $f, g \in \mathsf{H}^2(\mathsf{T}^n)$ we have

$$\langle \Delta f, g \rangle = \int \|k\|^2 \, f(k) \overline{\hat{g}}(k) dk = \int \hat{f}(k) \|k\|^2 \, \overline{\hat{g}}(k) dk = \langle f, \Delta g \rangle.$$

Immediately, the next result follows.

**Theorem 3.3** *The Schrödinger operator* $H = \Delta + V$, $V \in C^\infty(\mathsf{T}^n, \mathsf{R})$ *is self-adjoint on* $\mathsf{H}^2(\mathsf{T}^n)$ *where*

$$\mathsf{H}^2(\mathsf{T}^n) = \{f \in L^2(\mathsf{T}^n) : (1 + |k|^2)\hat{f}(k) \in l^2(\mathsf{Z}^n)\};$$

with

$$\hat{f}(k) = \frac{1}{(2\pi)^{\frac{n}{2}}} \int_{\mathsf{T}^n} f(x) e^{-ik.x} dx.$$

*Proof.* On Fourier space, $H = \Delta + V$ is a multiplication operator in $\mathsf{H}^2(\mathsf{T}^n)$ and thus self-adjoint there; c.f: [8].

Note, the spectrum $\sigma(\Delta) = \sigma_{ess}(H) = [0, \infty)$; [8].

Next, we review the Rayleigh-Schrödinger Perturbation Theory which is a major tool used in this paper.

# 4 The Rayleigh-Schrödinger Perturbation Theory

Consider a perturbed operator $T$ in the parameter $\varepsilon$ given by

$$T = T_0 + \varepsilon V \tag{15}$$

where $0 < \varepsilon \prec 1$, $T_0$ is a self-adjoint unperturbed operator and $V$ is a perturbation operator.

Suppose $T_0$ of (15) is an operator on a finite dimensional Hilbert space and has eigenvalue $\lambda_0$. $\lambda_0$ is called degenerate when the secular equation for $T_0$, that is, $\det(T_0 - \lambda) = 0$ has multiple roots at $\lambda_0$.

Following [4], finding eigenvalues of (15) is equivalent to solving the secular equation

$$det(T(\varepsilon) - \lambda) = (-1)[\lambda^n + a_1(\varepsilon)\lambda^{n-1} + \cdots + a_n(\varepsilon)] = 0$$

of degree $n$; and that if

$$F(\varepsilon, \lambda) := \lambda^n + a_1(\varepsilon)\lambda^{n-1} + \cdots + a_n(\varepsilon),$$

we have conditions for $\lambda = \lambda_0$ to be simple or degenerate. It specifies the Puiseux series in $(\varepsilon - \varepsilon_0)^{\frac{1}{p}}$ where $p$ is some positive integer as

$$\lambda_i^{(\varepsilon)} = \lambda_0 + \sum_{j=1}^{\infty} \alpha_j^i (\varepsilon - \varepsilon_0)^{\frac{j}{p_i}} \tag{16}$$

for some $m$ roots of $\lambda$ near $\lambda_0$, [4].

The Rayleigh-Schrödinger series for the eigenvalue $\lambda(\varepsilon)$ of the general case of linear operator (15) is given by

$$\lambda(\varepsilon) = \lambda_0 + \varepsilon \frac{\sum_{n=0}^{\infty} a_n \varepsilon^n}{\sum_{n=0}^{\infty} b_n \varepsilon^n} \tag{17}$$

with

$$a_n = \frac{(-1)^{n+1}}{2\pi i} \oint_{|\lambda - \lambda_0| = \alpha} (\Omega_0, [V(T_0 - \lambda)^{-1}]^{n-1} \Omega_0) \mathrm{d}\lambda$$

and

$$b_n = \frac{(-1)^{n+1}}{2\pi i} \oint_{|\lambda - \lambda_0| = \alpha} (\Omega_0, (T_0 - \lambda)^{-1})[V(T_0 - \lambda)^{-1}]^n \Omega_0 \mathrm{d}\lambda;$$

see e.g. [26,5] and [4].

The Rayleigh-Schrödinger perturbation procedure produces approximation to the eigenvalues and eigenvectors of the operator $T$ by sequence of successively higher order corrections to the eigenvalues and eigenvectors of the unperturbed operator $T_0$ knowing only those of $T_0$; for details, one may wish to see [27] and [28].

Consider now the eigenvalue problem

$$Tx_i = \lambda x_i \tag{18}$$

with $\lambda \neq 0$. Suppose the unperturbed eigenvalues $\lambda_i^0$, $(i = 1, 2, \cdots, n)$ are all distinct. Under this assumption, we have

$$\lambda_i(\varepsilon) = \sum_{k=0}^{\infty} \varepsilon^k \lambda_i^{(k)} \tag{19}$$

and

$$x_i(\varepsilon) = \sum_{k=0}^{\infty} \varepsilon^k x_i^{(k)}, i = 1, 2, \cdots, n \tag{20}$$

for sufficiently small $\varepsilon$. The zero-order terms $(\lambda_i^{(0)}, x_i^{(0)})$ are the eigenpair of the unperturbed operator $T_0$; that is

$$(T_0 - \lambda_i^{(0)} I) x_i^{(0)} = 0 \tag{21}$$

while it is assumed that $x_i^{(0)}$ which are mutually orthogonal have been normalized to unity so that

$$\lambda_i^{(0)} = \langle x_i^{(0)}, T_0 x_i^{(0)} \rangle.$$

Substituting (19) and (20) into (18) yields

$$(T_0 - \lambda_i^{(0)} I) x_i^{(k)} = -(V - \lambda_i^{(1)} I) x_i^{(k-1)} + \sum_{j=0}^{k-2} \lambda_i^{(k-j)} x_i^{(j)}, (j, k = 1, 2, \cdots, \infty), \tag{22}$$

$(i = 1, 2, \cdots, n)$. For fixed $i$ solvability of (22) requires that its right-hand-side be orthogonal to $x_i^{(0)}$ for all $k$; [27]. Thus, the value of $x_i^{(1)}$ is determined by $\lambda_i^{(j+1)}$. Specifically,

$$\lambda_i^{(j+1)} = \langle x_i^{(0)}, V x_i^{(j)} \rangle \tag{23}$$

where we have used the orthogonalization of $x_i^{(0)}$. The expession (23) is the required eigenvalue corrections.

Further eigenvector corrections $x_i^j$ are determined through $\lambda_i^{(2j+1)}$. For odd $j = 2j + 1$, $j = 0, 1, \cdots$, we have

$$\lambda_i^{(2j+1)} = \langle x_i^{(j)}, V x_i^{(j)} \rangle - \sum_{\mu=0}^{j} \sum_{\nu=1}^{j} \lambda_i^{(2j+1-\mu-\nu)} \langle x_i^{(\nu)}, x_i^{(\mu)} \rangle. \tag{24}$$

while for even $k = 2j, \, j = 1, 2, \cdots,$ we have

$$\lambda_i^{(2j)} = \langle x_i^{(j-1)}, V x_i^{(j)} \rangle - \sum_{\mu=0}^{j} \sum_{\nu=1}^{j} \lambda_i^{(2j-\mu-\nu)} \langle x_i^{(\nu)}, x_i^{(\mu)} \rangle. \tag{25}$$

The pair (24 , 25) is known as Dalgarno-Stewart Identity; [27].

Degenerate case arises when $T_0$ possesses multiple eigenvalues. In this situation, the straight-forward analysis presented above encounters series of complications. This is a consequence of the fact that the Rellich's series, (20), guarantees the existence of perturbation eigenvector expansion, (23), for certain special unperturbed eigenvectors only ; c.f: [4]. These special unperturbed eigenvectors cannot be specified a priori but must instead emerge from the perturbation procedure itself; [4].

Suppose

$$\lambda_1^{(0)} = \lambda_2^{(0)} = \cdots = \lambda_n^{(0)} = \lambda^0$$

of multiplicity $m$ corresponding to known orthonormal eigenvectors $x_1^{(0)}, x_2^{(0)}, \cdots, x_m^{(0)},$ and we are required to determine the appropriate linear combinations

$$y_i^{(0)} = a_1^{(i)} x_1^{(0)} + a_2^{(i)} x_2^{(0)} + \cdots + a_m^{(i)} x_m^{(0)}; (i = 1, \cdots, m) \tag{26}$$

so that the expansions of (19) and (20) are valid with $x_i^{(k)}$ replaced with $y_i^{(k)}$. Since we desire $\{y_i^{(0)}\}_{i=1}^{m}$ be orthonormal, then

$$a_1^{\mu} a_1^{\nu} + a_2^{\mu} a_2^{\nu} + \cdots + a_m^{\mu} a_m^{\nu} = \delta_{\mu,\nu} \tag{27}$$

where $\delta$ is the Kronecker delta. We can go ahead to find where this degeneracy is resolved.

Assume degeneracy is fully resolved at first order. That is here $\lambda_i^{(1)}, (i = 1, 2, \cdots, m)$, are all distinct. Then we determine $\{\lambda_i^{(1)}, y_i^{(0)}\}_{i=1}^{m}$ so that (22) be solvable for $k = 1$ and $i = 1, 2, \cdots, m$. In order for this to happen, it is both necessary and sufficient that for each fixed $i$,

$$\langle x_{\mu}^{(0)}, (V - \lambda_i^{(1)} I) y_i^{(0)} \rangle = 0; \mu = 1, 2, \cdots, m. \tag{28}$$

Substituting (26) into (28) and using orthonormality of $\{x_{\mu}^{(0)}\}_{\mu=1}^{m}$, we arrive at the following system of equations in matrix form [27]:

$$\begin{pmatrix} x_1^{(0)}, Vx_1^{(0)} \rangle & \dots & \langle x_1^{(0)}, Vx_m^{(0)} \rangle \\ \langle x_2^{(0)}, Vx_1^{(0)} \rangle & \dots & \langle x_2^{(0)}, Vx_m^{(0)} \rangle \\ \vdots & \ddots & \vdots \\ \langle x_m^{(0)}, Vx_1^{(0)} \rangle & \dots & \langle x_m^{(0)}, Vx_m^{(0)} \rangle \end{pmatrix} \begin{pmatrix} a_1^{(i)} \\ a_2^{(i)} \\ \vdots \\ a_m^{(i)} \end{pmatrix} = \lambda_i^{(i)} \cdot \begin{pmatrix} a_1^{(i)} \\ a_2^{(i)} \\ \vdots \\ a_m^{(i)} \end{pmatrix} \tag{29}$$

Thus, each $\lambda_i^{(1)}$ is an eigenvalue with corresponding eigenvector $(a_1^{(i)}, a_2^{(i)} \cdots, a_m^{(i)})^T$ of the matrix equation defined by $A_{\mu,\nu} = \langle x_\mu^{(0)}, A^{(1)}x_\nu^{(0)} \rangle$ and $(\mu, \nu = 1, 2, \cdots, m)$. By assumption, the symmetric matrix

$$\mathbf{A} = \begin{pmatrix} x_1^{(0)}, Vx_1^{(0)} \rangle & \dots & \langle x_1^{(0)}, Vx_m^{(0)} \rangle \\ \langle x_2^{(0)}, Vx_1^{(0)} \rangle & \dots & \langle x_2^{(0)}, Vx_m^{(0)} \rangle \\ \vdots & \ddots & \vdots \\ \langle x_m^{(0)}, Vx_1^{(0)} \rangle & \dots & \langle x_m^{(0)}, Vx_m^{(0)} \rangle \end{pmatrix}$$

has $m$ distinct eigenvalues and hence $m$ orthonormal eigenvectors described by (27).

Now that $\{y_i^{(0)}\}$ are fully determined, we have by equation (23) the following identities

$$\lambda_i^{(1)} = \langle y_i^{(0)}, Vy_i^{(0)} \rangle; (i = 1, 2, \cdots, m) \tag{30}$$

Furthermore, combining (27) and (29) gives

$$\langle y_i^{(0)}, Vy_j^{(0)} \rangle = 0; i \neq j. \tag{31}$$

The remaining eigenvalue corrections $\lambda_i^{(k)}; k \geq 2$ may be obtained similarly using the Dalgarno-Stewart Identity (24, 25).

Whenever (22) is solvable, we express its solution as

$$y_i^{(k)} = \hat{y}_i^{(k)} + \beta_{1,k}^{(i)} y_1^{(0)} + \beta_{2,k}^{(i)} y_2^{(0)} + \cdots + \beta_{m,k}^{(i)} y_m^{(0)}; \tag{32}$$

$(i = 1, \cdots, m)$; [27] and [28]. By the normalization, $\beta_{i,k}^{(i)} = 0; (i = 1, 2, \cdots, m)$ and $\beta_{j,k}^{(i)}$ are to be determined. The $\beta_{j,k}^{(i)}$ are given by

$$\beta_{j,k}^{(i)} = \frac{\langle y_j^{(0)}, Vy_j^{(0)} \rangle - \sum_{l=1}^{k-1} \lambda_i^{(k-l+1)} \beta_{j,l}^{(i)}}{\lambda_i^{(1)} - \lambda_j^{(1)}}; \tag{33}$$

for $i \neq j$ and zero otherwise. Higher order corrections are obtained similarly; see e.g. [27].

# 5 Simplicity of Spectrum of the Laplacian on the $n$-torus

High symmetry of $\mathsf{T}^n$ leads to different degrees of multiplicities of the Laplacian on smooth functions on $\mathsf{T}^n$; [28]. Here, "Multiplicity" means the number of ways a given number can be represented as sum of squares of $n$ integers. The multiplicity of an eigenvalue $\lambda$ of the Laplacian on $\mathsf{T}^n$ can easily be obtained with *Mathematica* using the command $\mathrm{SquaresR[n,\lambda]}$. The list of the possibilities is displayed with the command $\mathrm{PowersRepresentations[\lambda,n,n]}$. For example,

SquaresR[2,325] = 24  and

PowersRepresentations[325,2,2] = {1,18},{6,17},{10,15},

SquaresR[2,13] = 8,  and

PowersRepresentations[13,2,] = {0,13},{1,12},{2,11},{3,10},{4,9},{5,8},{6,7}.

We observe from equation ( 8) that on $\mathsf{T}^2$,

$$\sigma(\Delta) = k_1^2 + k_2^2 \tag{34}$$

with the corresponding normalized eigenfunction $\dfrac{1}{2\pi} e^{ik\cdot x}$ where $k = (k_1, k_2) \in \mathsf{Z}^2$ and $x \in \mathsf{R}^2$. The spectrum is therefore the set of the eigenvalues listed with their multiplicities $m$ thus:

$$\{(0, m=1), (1, m=4), (2, m=4), (4, m=4), (5, m=8), (8, m=4),$$

$$(9, m=4), \cdots, (13, m=8), \cdots, (25, m=12), \cdots, (125, m=16), \cdots\}.$$

Similarly for the 3-torus, we have

$$\sigma(\Delta) = k_1^2 + k_2^2 + k_3^2 \tag{35}$$

with normalized eigenfunction $\dfrac{1}{(2\pi)^{\frac{3}{2}}} e^{ik\cdot x}$. The spectrum is the set

$$\{(0, m=1), (1, m=6), (2, m=12), (3, m=8), (4, m=6), (5, m=24),$$

$$\cdots, (9, m=30), \cdots, (100, m=30), \cdots, (1000, m=144), \cdots\}.$$

Moreover, on $\mathsf{T}^4$, we have

$$\sigma(\Delta) = k_1^2 + k_2^2 + k_3^2 + k_4^2 \tag{36}$$

and

$$\sigma(\Delta) = \{(0, m = 1), (1, m = 8), (2, m = 24), (3, m = 32), (4, m = 24),$$

$$(5, m = 48), (6, m = 96), \cdots, (200, m = 744), \cdots, (2000, m = 3744), \cdots\}.$$

Therefore, for the n-torus, we have

$$\sigma(\Delta) = k_1^2 + k_2^2 + k_3^2 + k_4^2 + \cdots + k_n^2 \tag{37}$$

c.f: [1,9,10] and [12].

Now in what follows, we attempt to split the spectrum with appropriate potential. That is, with the potential $V$, we use the known eigenfunctions $x_i^{(k)} \in H^2(\mathsf{T}^n)$ based on the Dalgarno-Stewart identity to construct new eigenfunctions $y_i^{(0)} = a_m^{(i)} x_m^{(0)}$, $i = 1, 2, \cdots, m$ to replace $x_i^{(k)}$ such that $\sigma(H) = \lambda^{(j+1)}$ split. We must choose $x_m^{(0)}$ such that $a_m^\mu a_m^\nu = \delta_{\mu, \nu}$.

The main result of this paper is the following theorem.

**Theorem 5.1.** *There exists a* $2\pi$ *-periodic perturbation potential* $V \in C^\infty(\mathsf{T}^n, \mathsf{R})$ *such that the spectrum,* $\sigma(\Delta + V)$, *on the n-torus splits at first order.*

*Proof.* Let $V$ be a self-adjoint perturbation operator satisfying the assumptions (1) and (2) of theorem (5.1), then it follows that $\sigma(\Delta + V) \subset \mathsf{R}$. Let $\lambda$ be a non-zero eigenvalue of $\Delta$ on $\mathsf{T}^n$ with multiplicity $m$. By the Stewart-Dalgarno identity (24 , 25), choose orthonormal basis $x_m^{(0)}$ such that

$$\langle x_m^{(0)}, V x_n^{(0)} \rangle = \langle V x_m^{(0)}, x_n^{(0)} \rangle = \delta_{m,n}.$$

Then the first order correction matrix becomes

$$\begin{pmatrix} \langle x_1^{(0)}, V x_1^{(0)} \rangle & \cdots & \langle x_1^{(0)}, V x_m^{(0)} \rangle \\ \langle x_2^{(0)}, V x_1^{(0)} \rangle & \cdots & \langle x_2^{(0)}, V x_m^{(0)} \rangle \\ \vdots & \ddots & \vdots \\ \langle x_m^{(0)}, V x_1^{(0)} \rangle & \cdots & \langle x_m^{(0)}, V x_m^{(0)} \rangle \end{pmatrix} \begin{pmatrix} a_1^{(i)} \\ a_2^{(i)} \\ \vdots \\ a_m^{(i)} \end{pmatrix} = \lambda_i^{(i)} \begin{pmatrix} a_1^{(i)} \\ a_2^{(i)} \\ \vdots \\ a_m^{(i)} \end{pmatrix}$$

by the Rayleigh-Schrödinger perturbation theory.

Now, suppose degeneracy is fully resolved at first order. That is here $\lambda_i^{(1)}; (i = 1, 2, \cdots, m)$ are all distinct. Then we determine $\{\lambda_i^{(1)}, y_i^{(0)}\}_{i=1}^m$ so that (22) be solvable for $k = 1$ and $i = 1, 2, \cdots, m$. In order for this to happen, it is both necessary and sufficient that for each fixed $i$,

$$\langle x_\mu^{(0)}, (V - \lambda_i^{(1)} I) y_i^{(0)} \rangle = 0; \mu = 1, 2, \cdots, m. \tag{38}$$

Substituting (26) into (38) and using orthonormality of $\{x_\mu^{(0)}\}_{\mu=1}^m$, we arrive at the following system of equations in matrix form

$$\begin{pmatrix} \langle x_1^{(0)}, Vx_1^{(0)} \rangle & \cdots & \langle x_1^{(0)}, Vx_m^{(0)} \rangle \\ \langle x_2^{(0)}, Vx_1^{(0)} \rangle & \cdots & \langle x_2^{(0)}, Vx_m^{(0)} \rangle \\ \vdots & \ddots & \vdots \\ \langle x_m^{(0)}, Vx_1^{(0)} \rangle & \cdots & \langle x_m^{(0)}, Vx_m^{(0)} \rangle \end{pmatrix} \begin{pmatrix} a_1^{(i)} \\ a_2^{(i)} \\ \vdots \\ a_m^{(i)} \end{pmatrix} = A_\phi \begin{pmatrix} a_1^{(i)} \\ a_2^{(i)} \\ \vdots \\ a_m^{(i)} \end{pmatrix}$$

since $\mathsf{T}^n$ is fully determined by a matrix Lie group $A_\phi$ generated by $n$-tuple of distinct points $(\psi_1, \psi_2, \cdots, \psi_n) \in \mathsf{T}^n$.

Hence,

$$A_\phi \begin{pmatrix} a_1^{(i)} \\ a_2^{(i)} \\ \vdots \\ a_m^{(i)} \end{pmatrix} = \begin{pmatrix} \mathrm{e}^{i\phi_1} & 0 & \cdots & 0 \\ 0 & \mathrm{e}^{i\phi_2} & \cdots & 0 \\ \vdots & \vdots & \vdots & \vdots \\ 0 & \cdots & \cdots & \mathrm{e}^{i\phi_n} \end{pmatrix} \begin{pmatrix} a_1^{(i)} \\ a_2^{(i)} \\ \vdots \\ a_m^{(i)} \end{pmatrix} = \lambda_i^{(i)} \begin{pmatrix} a_1^{(i)} \\ a_2^{(i)} \\ \vdots \\ a_m^{(i)} \end{pmatrix}$$

Thus, each $\lambda_i^{(1)}$ is an eigenvalue with corresponding eigenvector $(a_1^{(i)}, a_2^{(i)} \cdots, a_m^{(i)})^T$ of the matrix equation $M$ defined by $A_{\mu,\nu} = \langle x_\mu^{(0)}, A^{(1)} x_\nu^{(0)} \rangle$, $(\mu, \nu = 1, 2, \cdots, m)$ and $A^1 := V$.

Since the choice of the distinct points are arbitrary and $\phi_1 \neq \phi_2 \neq \cdots \neq \phi_n$, the symmetric matrix

$$\mathbf{A}_\phi = \begin{pmatrix} \mathrm{e}^{i\phi_1} & 0 & \cdots & 0 \\ 0 & \mathrm{e}^{i\phi_2} & \cdots & 0 \\ \vdots & \vdots & \vdots & \vdots \\ 0 & \cdots & \cdots & \mathrm{e}^{i\phi_n} \end{pmatrix}$$

has $m$ distinct eigenvalues and hence $m$ orthonormal eigenvectors described by (27)

# 6 Illustration

We illustrate the main result with the following example. Consider the potential

$$V(x) = \sum_{t \in \mathsf{Z}^n} e^{-t_\alpha^2} e^{it.x} - 1 \tag{39}$$

where now, we understand

$$V(x) = 2 \sum_{t_1 \in \mathsf{Z}} \sum_{t_2 \in \mathsf{Z}} \cdots \sum_{t_n \in \mathsf{Z}} e^{-\|t\|_\alpha^2} \cos t.x. \tag{40}$$

Whereas $\alpha$ is a sequence of positive numbers, $t = k^2 + l^2$. We choose $\alpha$ such that the symmetry in the eigenvalue is resolved.

## 6.1 The unit circle

Let $\alpha = 1$ then

$$V(x) = \sum_{t \in \mathbb{Z}} e^{-t^2} e^{it.x} - 1. \tag{41}$$

Obviously, $V$ is real-valued and $2\pi$-periodic on $S^1$. On expansion for $t \in \mathbb{Z}$, we have

$$V(x) = 2\sum_{t=1}^{\infty} e^{-t^2} \cos tx.$$

Now following the Rayleigh-Schrödinger perturbation theory, we choose orthonormal basis $e_k$ and $e_l$ such that

$$\langle e_k, V e_l \rangle = \int_0^{2\pi} V e_{k-l} dx = \frac{1}{\sqrt{2\pi}} \int_0^{2\pi} V(x) e_{k-l} dx. \tag{42}$$

Let us consider the eigenvalue $\lambda = 1$ of the Laplacian $\Delta$ on $S^1$ which has multiplicity $m = 2$. We observe that the normalised eigenfunctions for the eigenvalues are $e_1 = \frac{1}{\sqrt{2\pi}} e^{ix}$ and $e_2 = \frac{1}{\sqrt{2\pi}} e^{-ix}$.

From (42) we have

$$a_{11} := \langle e_1, V(x) e_1 \rangle = 0$$

and

$$a_{12} := \langle e_1, V(x) e_2 \rangle = \int_0^{2\pi} \overline{e_1} V(x) e_2 dx$$

$$= \frac{2}{2\pi} \int_0^{2\pi} e^{-2ix} \sum_{t=0}^{\infty} e^{-t^2} \cos tx \, dx$$

$$= \cdots + 0 + 0 + \frac{e^{-4}}{\pi} \int_0^{2\pi} \cos^2(2x) dx = 2e^{-4}.$$

This leads us to form the required matrix, (on scaling $2e^{-4}$ by $\frac{1}{2}$ ),

$$\mathbf{A} = \begin{pmatrix} 0 & e^{-4} \\ e^{-4} & 0 \end{pmatrix}.$$

Clearly, $A$ has distinct eigenvalues $\mu_1 \approx -0.0183156$ and $\mu_2 \approx 0.0183156$.

Again, consider the eigenvalue $\lambda = 9$ with $m = 2$. Following the same procedure as the case of $\lambda = 1$ above, using its normalised eigenfunctions $e_1 = \dfrac{1}{\sqrt{2\pi}} e^{3ix}$ and $e_2 = \dfrac{1}{\sqrt{2\pi}} e^{-3ix}$, we arrive at the required matrix

$$\mathbf{B} = \begin{pmatrix} 0 & e^{-6} \\ e^{-6} & 0 \end{pmatrix}$$

which splits into distinct eigenvalues $\mu_1 \approx -0.00247875$ and $\mu_2 \approx 0.00247875$.

## 6.2 The 2-torus

Consider the eigenvalue $\lambda = 1$ of $\Delta$ on the 2-torus $\mathbf{T}^2 = S^1 \times S^1$ which has multiplicity $m = 4$. The normalised eigenfunctions for this eigenvalue are $e_1 = \dfrac{1}{2\pi} e^{-ix_1}$, $e_2 = \dfrac{1}{2\pi} e^{-ix_2}$, $e_3 = \dfrac{1}{2\pi} e^{ix_2}$ and $e_4 = \dfrac{1}{2\pi} e^{ix_1}$. The potential now becomes

$$V(x) = 2 \sum_{t \in \mathbf{Z}^2} e^{-\|t\|_\alpha^2} \cos tx \tag{43}$$

where here,

$$\langle e_k, V(x)e_l \rangle := \int_0^{2\pi} \int_0^{2\pi} \overline{e^{-i(k_1 x_1 + k_2 x_2)}} V(x) e^{-i(l_1 x_1 + l_2 x_2)} \, dx_1 dx_2. \tag{44}$$

Let $\alpha_1 = 1$ and $\alpha_2 = 2$. From equation (44) we obtain the matrix

$$\mathbf{C} = \begin{pmatrix} 1 & e^{-3} & e^{-3} & e^{-4} \\ e^{-3} & 1 & e^{-8} & e^{-3} \\ e^{-3} & e^{-8} & 1 & e^{-3} \\ e^{-4} & e^{-3} & e^{-3} & 1 \end{pmatrix}.$$

The matrix $C$ has all distinct eigenvalues $\mu_1 \approx 1.1093$, $\mu_2 \approx 0.999665$, $\mu_3 \approx 0.981684$, and $\mu_4 \approx 0.909346$.

Furthermore, consider the eigenvalue $\lambda = 5$ which has multiplicity $m = 8$. The normalised eigenfunctions are $e_1 = \dfrac{1}{2\pi} e^{-2ix_1 - ix_2}$, $e_2 = \dfrac{1}{2\pi} e^{-2ix_1 + ix_2}$, $e_3 = \dfrac{1}{2\pi} e^{-ix_1 - 2ix_2}$, $e_4 = \dfrac{1}{2\pi} e^{-ix_1 + 2ix_2}$, $e_5 = \dfrac{1}{2\pi} e^{ix_1 - 2ix_2}$, $e_6 = \dfrac{1}{2\pi} e^{ix_1 + 2ix_2}$, $e_7 = \dfrac{1}{2\pi} e^{2ix_1 - ix_2}$ and $e_8 = \dfrac{1}{2\pi} e^{2ix_1 + ix_2}$. From equation (44), we obtain the required matrix as

$$\mathbf{D} = \begin{pmatrix} 1 & e^{-8} & e^{-3} & e^{-19} & e^{-11} & e^{-27} & e^{-16} & e^{-24} \\ e^{-8} & 1 & e^{-19} & e^{-3} & e^{-27} & e^{-11} & e^{-24} & e^{-16} \\ e^{-3} & e^{-19} & 1 & e^{-32} & e^{-4} & e^{-36} & e^{-11} & e^{-27} \\ e^{-19} & e^{-3} & e^{-32} & 1 & e^{-36} & e^{-4} & e^{-27} & e^{-11} \\ e^{-11} & e^{-27} & e^{-4} & e^{-36} & 1 & e^{-32} & e^{-3} & e^{-19} \\ e^{-27} & e^{-11} & e^{-36} & e^{-4} & e^{-32} & 1 & e^{-19} & e^{-3} \\ e^{-16} & e^{-24} & e^{-11} & e^{-27} & e^{-3} & e^{-19} & 1 & e^{-8} \\ e^{-24} & e^{-16} & e^{-27} & e^{-11} & e^{-19} & e^{-3} & e^{-8} & 1 \end{pmatrix}.$$

$D$ has split eigenvalues $\mu_1 \approx 1.05993$, $\mu_2 \approx 1.05966$, $\mu_3 \approx 1.04165$, $\mu_4 \approx 1.04125$, $\mu_5 \approx 0.958717$, $\mu_6 \approx 0.958321$ $\mu_7 \approx 0.940371$ and $\mu_8 \approx 0.940099$

Furthermore, we consider the eigenvalue $\lambda = 125$ which has multiplicity $m = 16$. The normalised eigenfunctions are $e_1 = \frac{1}{2\pi} e^{-2ix_1 - ix_2}$, $e_2 = \frac{1}{2\pi} e^{-2ix_1 + ix_2}$, $e_3 = \frac{1}{2\pi} e^{-ix_1 - 2ix_2}$, $e_4 = \frac{1}{2\pi} e^{-ix_1 + 2ix_2}$, $e_5 = \frac{1}{2\pi} e^{ix_1 - 2ix_2}$, $e_6 = \frac{1}{2\pi} e^{ix_1 + 2ix_2}$, $e_7 = \frac{1}{2\pi} e^{2ix_1 - ix_2}$, $e_9 = \frac{1}{2\pi} e^{2ix_1 + ix_2}$. $e_9 = \frac{1}{2\pi} e^{-2ix_1 - ix_2}$, $e_{10} = \frac{1}{2\pi} e^{-2ix_1 + ix_2}$, $e_{11} = \frac{1}{2\pi} e^{-ix_1 - 2ix_2}$, $e_{12} = \frac{1}{2\pi} e^{-ix_1 + 2ix_2}$, $e_{13} = \frac{1}{2\pi} e^{ix_1 - 2ix_2}$, $e_{14} = \frac{1}{2\pi} e^{ix_1 + 2ix_2}$, $e_{15} = \frac{1}{2\pi} e^{2ix_1 - ix_2}$ and $e_{16} = \frac{1}{2\pi} e^{2ix_1 + ix_2}$.

Let $\alpha_1 = 0.01$ and $\alpha_2 = 0.1$. From equation 44, we obtain a $16 \times 16$ matrix as

$$\mathbf{E} = \begin{pmatrix} 1 & e^{-\frac{5}{8}} & \cdots & e^{-\frac{25}{121}} & e^{-\frac{25}{161}} \\ \vdots & \vdots & \ddots & \vdots & \vdots \\ e^{-\frac{25}{161}} & e^{-\frac{25}{121}} & \cdots & e^{-\frac{5}{8}} & 1 \end{pmatrix}.$$

$E$ has eigenvalues $\mu_1 \approx 3.01444$, $\mu_2 \approx 3.01441$, $\mu_3 \approx 1.54036$ $\mu_4 \approx 1.50697$, $\mu_5 \approx 1.32492$, $\mu_6 \approx 1.29193$, $\mu_7 \approx 0.770357$, $\mu_8 \approx 0.755073$, $\mu_9 \approx 0.690046$, $\mu_{10} \approx 0.651702$, $\mu_{11} \approx 0.498875$, $\mu_{12} \approx 0.481356$, $\mu_{13} \approx 0.204938$ $\mu_{14} \approx 0.204496$, $\mu_{15} \approx 0.0250688$ and $\mu_{16} \approx 0.0250487$.

### 6.3 The 3-torus

We demonstrate with the eigenvalue $\lambda = 1$ which on $T^3 = S^1 \times S^1 \times S^1$ has multiplicity $m = 6$. Its normalised eigenfunctions are $e_1 = \frac{1}{(2\pi)^{\frac{3}{2}}} e^{-ix_1}$, $e_2 = \frac{1}{(2\pi)^{\frac{3}{2}}} e^{-ix_2}$, $e_3 = \frac{1}{(2\pi)^{\frac{3}{2}}} e^{-ix_3}$, $e_4 = \frac{1}{(2\pi)^{\frac{3}{2}}} e^{ix_3}$, $e_5 = \frac{1}{(2\pi)^{\frac{3}{2}}} e^{ix_2}$ and $e_6 = \frac{1}{(2\pi)^{\frac{3}{2}}} e^{ix_1}$. The potential here is therefore

$$V(x) = 2 \sum_{t \in \mathbb{Z}^3} e^{-P_t P_\alpha^2} \cos tx. \tag{45}$$

Define

$$\langle e_k, V(x)e_l \rangle := \int_0^{2\pi} \int_0^{2\pi} \int_0^{2\pi} \overline{e^{-i(k_1 x_1 + k_2 x_2 + k_3 x_3)}} V(x) e^{-(l_1 x_1 + l_2 x_2 + l_3 x_3)} dx_1 dx_2 dx_3. \tag{46}$$

and take $\alpha_1 = 1$ and $\alpha_2 = 2$ as before, which leads to the matrix

$$\mathbf{F} = \begin{pmatrix} 1 & e^{-3} & e^{-1} & e^{-1} & e^{-3} & e^{-4} \\ e^{-3} & 1 & e^{-2} & e^{-2} & e^{-8} & e^{-3} \\ e^{-1} & e^{-2} & 1 & 1 & e^{-2} & e^{-1} \\ e^{-1} & e^9 & 1 & 1 & e^{-2} & e^{-1} \\ e^{-3} & e^{-8} & e^{-2} & e^{-2} & 1 & e^{-3} \\ e^{-4} & e^{-3} & e^{-1} & e^{-1} & e^{-3} & 1 \end{pmatrix}.$$

The eigenvalues of $F$ are $\mu_1 \approx 2.44993$, $\mu_2 \approx 0.999665$, $\mu_3 \approx 0.981684$, $\mu_4 \approx 0.948775$, $\mu_5 \approx 0.619943$ and $\mu_6 \approx 0$.

## 6.4 The 4-torus

Consider the eigenvalue $\lambda = 1$ of the $\Delta$ which has multiplicity $m = 8$ on $\mathsf{T}^4 = S^1 \times S^1 \times S^1 \times S^1$. It has the following normalised eigenfunctions $e_1 = \frac{1}{(2\pi)^2} e^{-ix_1}$, $e_2 = \frac{1}{(2\pi)^2} e^{-ix_2}$, $e_3 = \frac{1}{(2\pi)^2} e^{-ix_3}$, $e_4 = \frac{1}{(2\pi)^2} e^{-ix_4}$, $e_5 = \frac{1}{(2\pi)^2} e^{ix_4}$, $e_6 = \frac{1}{(2\pi)^2} e^{ix_3}$, $e_7 = \frac{1}{(2\pi)^2} e^{ix_2}$ and $e_8 = \frac{1}{(2\pi)^2} e^{ix_1}$.

The potential becomes

$$V(x) = 2 \sum_{t \in \mathbb{Z}^4} e^{-\|t\|_\alpha^2} \cos tx \tag{47}$$

where here, $x \in \mathbb{R}^4$ and

$$\langle e_k, V(x)e_l \rangle := \int_0^{2\pi} \int_0^{2\pi} \int_0^{2\pi} \int_0^{2\pi} \overline{e^{-i(k_1 x_1 + k_2 x_2 + k_3 x_3 + k_4 x_4)}} V(x) e^{-(l_1 x_1 + l_2 x_2 + l_3 x_3 + l_4 x_4)} dx_1 dx_2 dx_{3dx_4}. \tag{48}$$

The required matrix, when $\alpha_1 = 1$ and $\alpha_2 = 2$ is therefore

$$\mathbf{G} = \begin{pmatrix} 1 & e^{-3} & e^{-1} & e^{-1} & e^{-1} & e^{-1} & e^{-3} & e^{-4} \\ e^{-3} & 1 & e^{-2} & e^{-2} & e^{-2} & e^{-2} & e^{-8} & e^{-3} \\ e^{-1} & e^{-2} & 1 & 1 & 1 & 1 & e^{-2} & e^{-1} \\ e^{-1} & e^{-2} & 1 & 1 & 1 & 1 & e^{-2} & e^{-1} \\ e^{-1} & e^{-2} & 1 & 1 & 1 & 1 & e^{-2} & e^{-1} \\ e^{-1} & e^{-2} & 1 & 1 & 1 & 1 & e^{-2} & e^{-1} \\ e^{-3} & e^{-8} & e^{-2} & e^{-2} & e^{-2} & e^{-2} & 1 & e^{-3} \\ e^{-4} & e^{-3} & e^{-1} & e^{-1} & e^{-1} & e^{-1} & e^{-3} & 1 \end{pmatrix}.$$

The matrix $G$ has eigenvalues $\mu_1 \approx 4.37347$, $\mu_2 \approx 0.999665$, $\mu_3 \approx 0.981684$, $\mu_4 \approx 0.955542$, $\mu_5 \approx 0.689642$, $\mu_6 \approx -2.54159 \times 10^{-16}$, $\mu_7 \approx -5.67363 \times 10^{-17}$ and $\mu_8 \approx -4.2159 \times 10^{-17}$. Continuing this way, we see that the spectrum of the Laplace operator perturbed by this potential split at first order.

# 7 Conclusion

We proved the existence of a perturbation potential $V$ which guarantees the simplicity of the spectrum of the Schrödinger-type operator $\Delta + V$ on the n-torus. The Rayleigh-Schrödinger procedure of perturbation theory which basically produces approximation to the eigenvalues and eigenvectors of a perturbed operator by a sequence of successively higher order corrections to those of the unperturbed operator was employed to demonstrate that the spectrum of $\Delta + V$ splits generically at first order. With the properties prescribed on $V$, self-adjointness of the Laplacian carries over to $\Delta + V$.

The result is illustrated on different dimensional unit tori using the potential $V(x) = \sum_{t \in \mathbb{Z}^n} e^{-t_\alpha^2} e^{it.x} - 1$ which satisfies all the properties we require of $V$ on $\mathsf{T}^n$. We hope that the procedures of this paper can be followed to study generic spectrum simplicity of various perturbed Laplacians on other Riemannian manifolds.

# Competing Interests

Authors have declared that no competing interests exist.

# References

[1]  Grigor'yan A. Heat kernel and analysis on manifolds. American Mathematical Society AMS/PS Studies in advanced mathematics, USA; 2009.

[2]  Kiyosi *Itô*. Encyclopedic dictionary of mathematics. MIT Press, 5[th] Ed; 2000.

[3]  Courant R, Hilbert D. Methods of mathematical physics. John Wiley and Sons Inc; 1989.

[4]  Reed M, Simon B. Method of modern mathematical physics 4. Academic Press; 1980.

[5]     Kato T. Perturbation theory for Linear operators. Springer-Verlag, Berlin, Heidelberg; 1980.

[6]     Lions JL. Exact controllability, stabilization and perturbation for distributed systems. SIAM Rev. 1988;30(1):1-68.

[7]     Berezin FA, Shubin MA. The Schrödinger equation. Kluwer Academic Publishers, Netherlands; 1991.

[8]     Hislop PD, Sigal IM. Introduction to spectral theory: With applications to Schrödinger operators. Springer-Verlag Berlin, Germany; 1996.

[9]     Buser P. Geometry and spectra of compact riemannian surfaces. Birkh *ä* user Boston; 1992.

[10]    Chavel I. Eigenvalues in riemannian geometry. Academic Press Inc, London; 1984.

[11]    Omenyi L. Vacuum energy of the Laplacian on the spheres. Asian Research Journal of Mathematics. 2016;1(5):1-14. Article no. ARJOM.30523.
        Available: www.sciencedomain.org

[12]    Lee JM. Introduction to smooth manifolds. Graduate Text in Mathematics, Springer-Verlag, New York Inc; 2003.

[13]    Rayleigh JWS, Lindsay RB. The theory of sound. Second edition. Dover Classics of Science and Mathematics; 1945.

[14]    Hadamard J. Memoire sur le probl'eme danalyse relatif a lequilibre des plaques elastiques encastrees. Ouvres de J. Hadamard. 1968;II:515-641.

[15]    Albert JH. Genericity of simple eigenvalues for elliptic PDEs. Proc. Amer. Math. Soc. 1975;48:413-418.

[16]    Micheletti AM. Perturbazione dello spectro dell operatore de Laplace in relazione ad una variazone del campo. Ann. Scuola Norm. Sup. Pisa Cl. Sci. 1972;26(3):151-169.

[17]    Uhlenbeck K. Generic properties of eigenfunctions. Amer. J. Math. 1976;98:1059-1078.

[18]    Henry DB. Perturbation of the boundary in boundary-value problems of partial differential equations. London Mathematical Society, Cambridge University Press, UK, NY, USA; 2005.

[19]    Lagnese JE. Boundary stabilization of thin plates. SIAM Stud. Appl. Math; 1989.

[20]    Ortega JH, Enrique Zuazua. Generic simplicity of the spectrum and stabilization for a plate equation. SIAM. J. Control Optim. 2000;39(15):1585-1614.

[21]    Marcone CP. Generic simplicity of eigenvalues for Dirichlet problem of the biLaplacian operator. Electronic Journal of Differential Equations. 2004;2004(114):1-21.

[22]    Hillairet Luc. Generic simplicity of polygons. Workshop on Spectrum and Dynamics. 2008;7-11.

[23]    Dobrokhotov SY, Shafarevich AI. "Momentum" tunneling between Tori and the splitting of eigenvalues of the Laplace-Beltrami operator on Louville surface. Mathematical Physics, Analysis and Geometry. 1999;2(2):141-177.

[24] Dabrowski A. Small volume expansion of the splitting of multiple Neumann Laplacian eigenvalues due to a grounded inclusion in two dimensions. arXiv:1612.09095v1; [Math.AP], 29 December; 2016.

[25] Sarnak P. Spectra of hyperbolic surfaces. Bull. Amer. Math. Soc. (N.S.). 2003;40(4):441-478.

[26] Millman RS. Remarks on spectrum of the Laplace-Beltrami operator in the middle dimensions. Tensor 34. 1980;94-96.

[27] Brian JM. Rayleigh-Schrödinger perturbation theory: Pseudoinverse formulation. Hikari Ltd; 2009.

[28] Parlett BN. The symmetric eigenvalue problem. Society for Industrial and Applied Mathematics; 1998.